\documentclass{article}
\usepackage{latexsym}
\usepackage{amsmath,amsthm,amssymb}
\newcommand{\brac}[1]{\left(#1\right)}

\def\cA{{\cal A}}

\def\E{{\bf E}}

\def\a{\alpha}

\def\g{\gamma}

\def\z{\zeta}

\def\l{\lambda}

\def\p{\pi}

\def\om{\omega}
\def\OM{\Omega}

\def\Pr{\mbox{{\bf P}}}

\def\whp{{\bf whp}}
\def\Whp{{\bf Whp}}

\newtheorem{lemma}{Lemma}
\newtheorem{theorem}{Theorem}

\newcommand{\proofstart}{{\bf Proof\hspace{2em}}}

\newcommand{\proofend}{\hfill\mbox{$\Box$}}

\def\a{\alpha}

\def\Pr{\mbox{{\bf Pr}}}

\begin{document}

\large

\makeatletter \title{The emergence of a giant component in random subgraphs
of pseudo-random graphs}
\author{Alan Frieze\thanks{
Department of Mathematical Sciences, Carnegie Mellon University,
Pittsburgh PA15213, U.S.A. Supported in part by NSF grant CCR-9818411.}
\and Michael Krivelevich\thanks{
Department of Mathematics, Raymond and Beverly Sackler Faculty of Exact
Sciences, Tel Aviv University, Tel Aviv 69978, Israel. E-mail: 
krivelev@math.tau.ac.il. Research supported in part by USA-Israel
BSF Grant 99-0013, by grant 64-01 from the Israel Science Foundation 
and by a Bergmann Memorial Grant.}
\and Ryan Martin\thanks{Department of Mathematical Sciences,
Carnegie Mellon University, Pittsburgh, PA 15213. Supported in part by NSF
VIGRE grant DMS-9819950.}}\date{}
\maketitle \makeatother
\begin{abstract}
Let $G$ be a $d$-regular graph $G$ on $n$ vertices. Suppose that 
the adjacency matrix of $G$ is such that the
eigenvalue $\l$ which is second largest in absolute value satisfies $\l=o(d)$. Let
$G_p$ with $p=\frac{\a}{d}$ be obtained from $G$ by including each edge of
$G$ independently with probability $p$. We show that if $\a<1$ then
\whp\ the maximum component size of $G_p$ is $O(\log n)$ and if $\a>1$
then $G_p$ contains a unique giant component of size $\OM(n)$, with
all other components of size $O(\log n)$.
\end{abstract}

\section{Introduction}
Pseudo-random graphs (sometimes also called quasi-random graphs) can be
informally defined as graphs whose edge distribution resembles closely
that of truly random graphs on the same number of vertices and with the
same edge density. Pseudo-random graphs, their constructions and
properties have been a subject of intensive study for the last fifteen
years (see \cite{AloKriSud99}, \cite{ChuGraWil89}, \cite{Tho87a}, 
\cite{SimSos97}, \cite{Tho87b}, to mention just a few).

For the purposes of this paper, a pseudo-random graph is a $d$-regular
graph $G=(V,E)$ with vertex set $V=[n]=\{1,\ldots,n\}$,
all of whose eigenvalues but the first one are
significantly smaller than $d$ in their absolute values. More formally,
let $A=A(G)$ be the adjacency matrix of $G$. This is an $n$-by-$n$
matrix such that $A_{ij}=1$ if $(i,j)\in E(G)$ and $A_{ij}=0$
otherwise. Then $A$ is a real symmetric matrix with non-negative values
of its entries. Let $\lambda_1\ge\lambda_2\geq \cdots\ge\lambda_n$ be the
eigenvalues of $A$, also called the eigenvalues of $G$.
It follows from the Perron-Frobenius theorem that $\lambda_1=d$ and
$|\lambda_i|\le d$ for all $2\le i\le n$. We thus denote
$\lambda=\lambda(G)=\max_{2\le i\le n} |\lambda_i|$. The reader is
referred to a monograph of Chung \cite{Chu97} for further information on
spectral graph theory.

It is known (see, e.g. \cite{AloChu88}) that the greater is the so-called
spectral gap (i.e. the difference between $d$ and $\lambda$) the more
tightly the distribution of the edges of $G$ approaches that of the
random graph $G(n,d/n)$. We will cite relevant quantitative results
later in the text (see Lemma \ref{concen}),
 for now we just state informally that a spectral gap ensures
pseudo-randomness. 

In this paper we study certain properties of a random subgraph of a
pseudo-random graph. Given a graph $G=(V,E)$ and an edge probability
$0\le p=p(n)\le 1$, the {\em random subgraph} $G_p$ is formed by
choosing
each edge of $G$ independently and with probability $p$. 
We will also need to consider the related random graph $G_m$ whose
edge set is a random $m$-subset of $E$.

The most
studied random graph is the so called binomial random graph $G(n,p)$,
formed by choosing the edges of the complete graph on $n$ labeled
vertices independently with probability $p$. Here rather than studying
random subgraphs of one particular graph, we investigate the properties
of random subgraphs of graphs from a wide class of regular pseudo-random
graphs. As we will see, all such subgraphs viewed as probability spaces
share certain common features. 

Our concern here is with the existence of a giant component in the
case $p=\frac{\a}{d}$ or $m=\frac{1}{2}\a n$ where
$\a\neq 1$
is an absolute constant. These two models are sufficiently similar
so that the results we prove in $G_p$ immediately translate to $G_m$
and vice-versa. The needed formal relations in the case where $G=K_n$
are given in \cite{Boll} or \cite{JanLucRuc00} and they generalise
easily to our case.

As customary when studying random graphs, asymptotic conventions and
notations apply. In particular, we assume where necessary the number of
vertices $n$ of the base graph $G$ to be as large as needed. Also, we
say that a
graph property $\cA$ holds {\em with high probability}, or \whp\ for
brevity, in $G_p$ if the probability that $G_p$ has $\cA$ tends to 1 as
$n$ tends to infinity. Monographs \cite{Boll}, \cite{JanLucRuc00} provide
a necessary background and reflect the state of affairs in the theory
of random graphs.

For $\a>1$ we define $\bar{\a}<1$ to be the unique solution (other than $\a$) 
of the equation $xe^{-x}=\a e^{-\a}$. 
We assume from now on that
\begin{equation}\label{bound}
d\to\infty\mbox{ and }\l=o(d).
\end{equation}
These requirements are quite minimal.

In analogy to the classical case $G=K_n$, studied already by Erd\H
os and R\'enyi \cite{ER1},
\begin{theorem}\label{th1}
Assume that \eqref{bound} holds.
\begin{itemize}
\item[(a)] If $\a<1$ then \whp\ the maximum component size is $O(\log n)$.  
\item[(b)] If $\a>1$ then \whp\ there is a unique giant component of
asymptotic size $\brac{1-\frac{\bar{\a}}{\a}}n$ and the remaining components
are of size $O(\log n)$.
\end{itemize}
\end{theorem}
One can also prove tighter results on the size and structure of the
small components. They correspond nicely to the case where $G=K_n$.

We will use the notation $f(n)\gg g(n)$ to mean $f(n)/g(n)\to\infty$
with $n$. Similarly, $f(n)\ll g(n)$ means that $f(n)/g(n)\to 0$.
\begin{theorem}\label{th2}
Assume that \eqref{bound} holds. Let $\om=\om(n)\to\infty$ with $n$.
\begin{itemize}
\item[(a)] If $d\gg(\log n)^2$ then {\bf whp} $G_p$ contains no isolated
 trees of size $\z(\log
n-\frac{5}{2}\log\log n)+\om$, where
$\z^{-1}=\a-1-\log \a>0$.
\item[(b)] If $d\gg \log^2 n$ then {\bf whp} $G_p$ contains an isolated tree of size 
at least $\z(\log n-\frac{5}{2}\log\log n)-\om$.
\item[(c)] If $d=\OM(n)$ then \whp\ 
$G_p$ contains $\leq \om$ vertices on unicyclic
components.
\item[(d)] Let $d\gg \sqrt{n}$.
If $\a<1$ then \whp\ $G_p$ contains no component with $k$
vertices and with more than $k$ edges.
\item[(e)] Let $d\gg \sqrt{n}$.
If $\a>1$ then \whp\ $G_p$ contains no component with $k=o(n)$
vertices and with more than $k$ edges.
\end{itemize}
\end{theorem}
\section{Properties of $d$-regular graphs}

In this section we put together those properties needed to prove
Theorem~\ref{th1}.  For $B, C\subseteq V$ 
let $e(B,C)$ denote the number of {\bf ordered} pairs
$(u,v)$ such that $u\in B$, $v\in C$ and $\{u,v\}\in E$.

\begin{lemma}\label{concen}
Suppose $B, C\subseteq V$ and $|B|=bn$ and $|C|=cn$. Then
\[ \left|e(B,C)-bcdn\right|\leq\lambda n\sqrt{bc}. \]
\end{lemma}

This is Corollary 9.2.8 of \cite{AS}.  Note that $B=C$ is allowed here.
Then $e(B,B)$ is twice the number of edges of $G$ in the graph induced by
$B$.

Now let $t_k$ denote the number of $k$-vertex trees that are contained
in $G$.

\begin{lemma}
\[ n\frac{k^{k-2}(d-k)^{k-1}}{k!}\leq t_k\leq n\frac{k^{k-2}d^{k-1}}{k!} \]
\label{tkbound}
\end{lemma}

This is Lemma 2 of \cite{BFM}.

\section{Proof of Theorem~\ref{th1}}

Let $p=\frac{\alpha}{d}$
and let $C_k$ denote the number of vertices of $V$ that are contained
in components of size $k$ in $G_p$ and let $T_k\leq C_k$ denote the
number of vertices which are contained in isolated trees of size $k$.

\begin{lemma}\label{lem1} \ 
\begin{itemize}
\item[(a)] 
$${\bf E}C_k\leq n\frac{k^{k-1}}{k!}\alpha^{k-1}e^{-\a k(1-\xi_k)}$$ 
where
$$\xi_k=\min\left\{\frac{k}{d}, \frac{
k}{n}+\frac{\lambda}{d}\right\}.$$
\item[(b)] For $k\ll d$,
$${\bf E}T_k\geq
n\frac{k^{k-1}}{k!}\alpha^{k-1}e^{-\a k(1+\eta_k)}$$ 
where
$$\eta_k=\frac{2k}{d}+\frac{2k}{\a d}+\frac{\a}{d}.$$
\end{itemize}
\label{explemma}
\end{lemma}

\proofstart
\begin{itemize}
\item[(a)] Let ${\cal T}_k$ denote the set of trees of size $k$ in $G$.  Then
\[ {\bf E}C_k\leq\sum_{T\in {\cal T}_k}kp^{k-1}(1-p)^{e_T} \]
where $e_T=e\left(V(T),\overline{V(T)}\right)$.  Now
Lemma~\ref{concen} implies that
\begin{equation}
   e_T=kd-e(V(T),V(T))\geq a_k\stackrel{\rm def}{=}
   kd-\frac{k^2d}{n}-\l k
   \label{etbound}
\end{equation}
and we also have the simple inequality
\[ e_T\geq b_k\stackrel{\rm def}{=}kd-k(k-1) \]
which is true for an arbitrary $d$-regular graph.

Thus,
\begin{equation}\label{3}
   {\bf E}C_k\leq kt_kp^{k-1}(1-p)^{\max\{a_k,b_k\}}
\end{equation}
and (a) follows from Lemma~\ref{tkbound} and some straightforward estimations.

\item[(b)] Similarly,
\[ {\bf E}T_k\geq\sum_{T\in {\cal T}_k}kp^{k-1}
   (1-p)^{e_T+k^2} \]
where we crudely bound by $k^2$, the number of 
edges contained in $V(T)$ which must be absent to make $T$ an isolated
tree component of $G_p$.  Now we can simply use
$$e_T\leq kd$$
and Lemma~\ref{tkbound}.  We also use $1-p\geq e^{-p-p^2}$ for $p$
small and
\[ (d-k)^{k-1}>d^{k-1}\left(1-\frac{k}{d}\right)^k\geq 
   d^{k-1}\exp\left\{-\frac{k^2}{d}-\frac{k^3}{d^2}\right\}\; , \]
for $k/d$ small and make some straightforward estimations. 
\end{itemize}
\proofend

Now choose 
$\g=\g(\a)$ such that 
$$\a e^{1-\a+2\a\g}=1.$$
(Note that $\alpha e^{1-\alpha}<1$ for $\alpha\neq 1$.)

\begin{lemma}\label{lem3}
  \Whp, $C_k=0$ for $k\in I=\left[\frac{1}{\a\g}\log n, \gamma 
  n\right]$
  \label{size}
\end{lemma}

\proofstart 
First assume that $k\leq \g d$ and observe that $\xi_k\leq \g$ in
this range. Then from Lemma \ref{lem1}(a) and 
$k!\geq\left(\frac{k}{e}\right)^k$ we see that
\begin{eqnarray}
\sum_{k=\frac{1}{\a\g}\log n}^{\g d}{\bf E}\,C_k&\leq& 
\frac{n}{\alpha}\sum_{k=\frac{1}{\a\g}\log n}^{\g d}
k^{-1}(\a e^{1-\alpha+\a\xi_k})^k\nonumber\\
&\leq&\frac{n}{\alpha}\sum_{k=\frac{1}{\a\g}\log n}^{\g d}
k^{-1}(\a e^{1-\alpha+\a\g})^k\nonumber\\
&=&\frac{n}{\alpha}\sum_{k=\frac{1}{\a\g}\log n}^{\g d}k^{-1}e^{-\a\g k}
\nonumber\\
&\leq&\frac{\g n}{\log n}\sum_{k=\frac{1}{\a\g}\log n}^{\infty}e^{-\a\g
k}\nonumber\\
&=&o(1).\label{f4}
\end{eqnarray}
Now assume that $\g d\leq k\leq \g n$ and observe that 
(\ref{bound}) implies $\xi_k\leq \g+o(1)$ in
this range. Then
\begin{eqnarray}
\sum_{k=\frac{1}{\a\g}\log n}^{\g n}{\bf E}\,C_k&\leq&
\frac{n}{\alpha}\sum_{k=\frac{1}{\a\g}\log n}^{\g n}k^{-1}(\a e^{1-\alpha+\a\g+o(1)})^k\nonumber\\ 
&\leq&\frac{n}{\alpha}\sum_{k=\frac{1}{\a\g}\log n}^{\g n}k^{-1}e^{-(\a\g-o(1)) k}\nonumber\\
&=&o(1).\label{f3}
\end{eqnarray}
\proofend

Now let us show that there are many vertices on small isolated trees.
Let
$$f(\a)=\sum_{k=1}^{\infty}\frac{k^{k-1}}{k!}\alpha^{k-1}e^{-\a k}.$$
It is known, see for example Erd\H{o}s and R\'enyi \cite{ER1} that
$$f(\a)=\begin{cases}1&\a\leq 1.\\ \frac{\bar{\a}}{\a}&\a>1.\end{cases}$$
\begin{lemma}\label{lem2}
Let $k_0=d^{1/3}$. Then
$$\Pr\brac{\left|\sum_{k=1}^{k_0} C_k-nf(\a)\right|\geq n^{5/6}\log
n}=o(1).$$
\end{lemma}
\proofstart
Note that $k\xi_k=O(d^{-1/3})$ and 
$k\eta_k=O(d^{-1/3})$ for $k\leq k_0$. Thus from Lemma \ref{lem1}(a) we have 
\begin{equation}\label{f2}
\E\sum_{k=1}^{k_0}C_k\leq
(1+O(d^{-1/3}))n\sum_{k=1}^{k_0}
\frac{k^{k-1}}{k!}\alpha^{k-1}e^{-\a k}=(1+O(d^{-1/3}))nf(\a).
\end{equation}
On the other hand, Lemma \ref{lem1}(b) implies,
\begin{equation}\label{f1}
\E\sum_{k=1}^{k_0}C_k\geq
\E\sum_{k=1}^{k_0}T_k\geq (1-O(d^{-1/3}))n\sum_{k=1}^{k_0}
\frac{k^{k-1}}{k!}\alpha^{k-1}e^{-\a k}=(1-O(d^{-1/3}))nf(\a).
\end{equation}
We now use the Azuma-Hoeffding martingale tail inequality \cite{AS} to 
show that the random variable $Z=\sum_{k=1}^{k_0}C_k$ is sharply
concentrated. We switch to the model $G_m,m=\frac{1}{2}\a n$. Changing
one edge can only change $Z$ by at most $2k_0$ and so for any $t>0$
$$\Pr(|Z-\E Z|\geq t)\leq 2\exp\left\{-\frac{2t^2}{4mk_0^2}\right\}.$$
Putting $t=n^{1/2}k_0\log n$ yields the lemma, in conjunction with
(\ref{f2}), (\ref{f1}) and $d\to\infty$.
\proofend

The first part of Theorem \ref{th1} now follows easily. Since $\a<1$
here, we have $f(\a)=1$ and so by Lemma \ref{lem3} and Lemma
\ref{lem2} \whp\ there are at least $n-n^{5/6}\log n$ vertices in
components of size at most $\frac{1}{\a\g}\log n$. Applying Lemma
\ref{lem3} again, we see that \whp\ the remaining vertices $X$ must be in
components of size at least $\g n$. So if $X\neq \emptyset$ then
$|X|\geq \g n$. But we know that \whp\ $|X|\leq n^{5/6}\log n$ and so
$X=\emptyset$ \whp.

For the second part of the theorem where $\a>1$ we see that \whp\
there are\\ $\frac{\bar{\a}}{\a}n+O(n^{5/6}\log n)$ vertices on
components of size $\leq\frac{1}{\a\g}\log n$ and the remaining
vertices lie in {\em large} components of size at least $\g n$. This statement
remains true if we consider $G_{m-\log n}$. Let $S_1,S_2,\ldots,S_s$
be the large components of $G_{m-\log n}$, where $s\leq 1/\g$. We now
show that \whp\ adding the remaining $\log n$ random edges $Y$ puts
$S_1,S_2,\ldots,S_s$ together in one giant component of size
$\brac{1-\frac{\bar{\a}}{\a}}n+O(n^{5/6}\log n)$. We also \whp\ have 
$\frac{\bar{\a}}{\a}n+O(n^{5/6}\log n)$ vertices on
components of size $\leq\frac{1}{\a\g}\log n$ and Lemma \ref{lem3}
shows that this accounts for all the vertices.

So let us show that
\begin{equation}\label{ryan}
\Pi=\Pr(\exists\, 1\leq i<j\leq s:\;Y\mbox{ contains no edge joining 
$S_i$ and $S_j$})=o(1),
\end{equation}
completing the proof of Theorem \ref{th1}.
Now by Lemma \ref{concen}, $G$ contains at least $(1-o(1))\g^2dn$
edges between $S_i$ and $S_j$, 
and the probability that $Y$ contains none of these is at most
$\brac{1-\frac{(1-o(1))\g^2dn}{\frac{1}{2}dn}}^{\log n}\leq
n^{-2\g^2+o(1)}$.
So $\Pi\leq \g^{-2}n^{-2\g^2+o(1)}=o(1)$, proving (\ref{ryan}).
\proofend
\section{Proof of Theorem \ref{th2}}
Let $k_\pm=\z(\log n-\frac{5}{2}\log\log n)\pm\om$. Let $N_k$ denote
the number of tree components of size $k$ in $G_p$.

Assume that $k_-\leq k\leq\frac{1}{\a\g}\log n$. Then from {\em the
proof of} Lemma \ref{lem1}, (notice $k^{k-2}$ in place of $k^{k-1}$,
we are counting trees, not vertices on trees),
\begin{eqnarray}
\E N_k&\leq&n\frac{k^{k-2}}{k!}\a^{k-1}e^{-\a k(1-\xi_k)}\nonumber\\
&=&(1+o(1))\frac{n}{\a k^{5/2}\sqrt{2\p}}e^{-\z^{-1}k}\label{uppera}\\
\noalign{and when $k=o(d^{1/2})$}
\E N_k&=&(1+o(1))\frac{n}{\a k^{5/2}\sqrt{2\p}}e^{-\z^{-1} k}\label{upper}\qquad
\end{eqnarray}
(a) Let $\g$ be as in Lemma \ref{lem3}. Using \eqref{uppera},
$$\sum_{k=k_+}^{\frac{1}{\a\g}\log n}\E
N_k=O\brac{\sum_{k=k_+}^{\frac{1}{\a\g}\log n}e^{-\z^{-1}(\om+k-k_+)}}=o(1)$$
and part (a) will follow once we verify that when $\a>1$, the giant
component is not a tree. However, the number of edges in the giant is
asymptotically 
\begin{eqnarray*}
\frac{\a n}{2}-\frac{n}{\a}\sum_{k=1}^\infty \frac{(k-1)k^{k-2}}{k!}(\a
e^{-\a})^k&=&\a n\brac{\frac{1}{2}-\frac{1}{\a^2}\sum_{k=1}^\infty
\frac{(k-1)k^{k-2}}{k!}(\bar{\a} e^{-\bar{\a}})^k }\\
&=&\a n\brac{\frac{1}{2}-\frac{\bar{\a}^2}{2\a^2}}.
\end{eqnarray*}
Note that 
$$\frac{n}{\bar{\a}}
\sum_{k=1}^\infty\frac{(k-1)k^{k-2}}{k!}(\bar{\a} e^{-\bar{\a}})^k
=\frac{\bar{\a}}{2}n$$
which can be seen from the fact that the LHS is asymptotically equal
to the expected number of edges of $G_{n,\frac{\bar{\a}}{n}}$ which
lie on trees. So, the ratio of edges to vertices for the giant is
asymptotically equal to $\frac{\a+\bar{a}}{2}>1$.

(b) Now let $k=k_-$. Then from (\ref{upper}),
$$\E N_k=\Omega(e^{\z^{-1}\om})\to\infty.$$
Bounding the number of $G$-edges
inside and between two disjoint subtrees by $3k^2$ we estimate
\begin{eqnarray*}
\E N_k^2&\leq&t_k^2p^{2k-2}(1-p)^{2dk-3k^2}\\
&=&(1+o(1))(\E N_k)^2
\end{eqnarray*}
and (b) follows from the Chebychev inequality. 

(c) Let $U_k$ denote the number of isolated unicyclic components
in $G_p$ of size $k$. Then
\begin{eqnarray*}
\sum_{k=3}^n\E (kU_k)&\leq&\sum_{k=3}^nkt_k\binom{k}{2}p^k(1-p)^{dk-k^2}\\
&\leq&(1+o(1))\frac{n}{2d}\sum_{k=3}^n\frac{k^{k+1}}{k!}(\a e^{-\a+o(1)})^k\\
&\leq&(1+o(1))\frac{n}{2d}\sum_{k=3}^n\frac{k^{1/2}}{\sqrt{2\p}}(\a e^{1-\a+o(1)})^k\\
&=&O(1)
\end{eqnarray*}
since we are assuming that $d=\OM(n)$ here. 
Part (c) follows from the Markov inequality.

(d), (e) Let $COMP_k$ denote the number of components with $k$ vertices and
at least $k+1$ edges. We can restrict our attention to $4\leq k\leq\g n$
since if $\a<1$ there are no larger components \whp. Then, as in
(\ref{3}),
\begin{eqnarray*}
\E \sum_{k=4}^{\g n} COMP_k&\leq& \sum_{k=4}^{\g n}
t_k\binom{k}{2}^2p^{k+1}(1-p)^{\max\{a_k,b_k\}}\\
&\leq&(1+o(1))\frac{n\a}{4\sqrt{2\pi}d^2}\sum_{k=4}^{\g n}k^{3/2}(\a
e^{1-\a+\a\g+o(1)})^k.\\
&\leq&(1+o(1))\frac{n\a}{4\sqrt{2\pi}d^2}\sum_{k=4}^{\g n}k^{3/2}
e^{-(\a\g-o(1)) k}\\
&=&o(1)
\end{eqnarray*}
since $n/d^2\to 0$. 

This completes the proof of (d), (e). 
\proofend

\end{document}